\documentclass[12pt]{article}

\usepackage{amsmath}
\usepackage{amssymb}
\usepackage{amsthm}

\newtheorem{theorem}{Theorem}

\newtheorem{claim}{Claim}

\newtheorem{conjecture}{Conjecture}

\newcommand{\lra}{\leftrightarrow}
\newcommand{\ra}{\rightarrow}
\newcommand{\Del}{\Delta}
\newcommand{\Z}{{\mathbb Z}}

\newcommand{\R}{{\mathbb R}}

\newcommand{\bv}{{\bf v}}
\newcommand{\bz}{{\bf z}}

\newcommand{\bw}{{\bf w}}
\newcommand{\be}{{\bf e}}
\newcommand{\bx}{{\bf x}}

\title{Simulating a Random Walk with Constant Error}
\author{JOSHUA N. COOPER
    \thanks{Research supported by NSF Grant DMS-0303272.}
    \@ and \@ JOEL SPENCER\\
    Department of Mathematics, Courant Institute of Mathematics\\
    New York University, \\
    251 Mercer St, New York, NY 10012--1185 \\
    cooper@cims.nyu.edu and spencer@cs.nyu.edu}

\begin{document}

\maketitle

\begin{abstract}
We analyze Jim Propp's $P$-machine, a simple deterministic process that simulates a random walk on $\Z^d$ to within a constant.  The proof of the error bound relies on several estimates in the theory of simple random walks and some careful summing.  We mention three intriguing conjectures concerning sign-changes and unimodality of functions in the linear span of $\{p(\cdot,\bx) : \bx \in \Z^d\}$, where $p(n,\bx)$ is the probability that a walk beginning from the origin arrives at $\bx$ at time $n$.
\end{abstract}

Consider the following machine, which we call the $P$-machine after its progenitor, Jim Propp.  
We begin with a set of chips on ``even'' vertices of $\Z^d$, i.e., those an even $L^1$-distance from the origin.  We label by 
$\Z_{2d}=\{0,1,\ldots,2d-1\}$
the $2d$ directions of $Z^d$.  (For example: North, South, East and West for $d=2$.)
On each vertex $\bv$ there will be a ``rotor" which will have a state $j\in \Z_{2d}$.
The initial states of the rotors can be set arbitrarily.  The rotors work as
follows:  
when one feeds it a chip, a rotor changes to state $j+1$ (addition
in $\Z_{2d}$) and sends the 
chip to the vertex in direction $j$.  Each time unit we feed all of the chips 
currently on each vertex to the rotor on that vertex.  The result is a new
distribution of chips.  We imagine doing this operation for total time $n$.

Now, consider the following process.  We again begin with a set of chips on ``even'' vertices of $\Z^d$.  Each chip takes an $n$-step simple random walk from its starting point.  At time $n$, we 
expect that the $P$-machine and the random process should give rise to similar distributions 
if they begin with the same configuration of chips.  For the random walk process we may
consider the expected number of chips that will be at $\bv$ at time $n$.  Our main result
is that the difference between this expected number and the actual number at $\bv$ at
time $n$ in the deterministic $P$-machine
is bounded uniformly -- {\it 
irrespective of how much time has passed, what the original chip distribution was, the 
starting states of the rotors, or even the choice of $\bv$!}

As an example, suppose $n$ is even, $d=1$, we begin with $n$ chips at position $0$, and the total
time is $n$.  The random walk model will have an expected number $n{n\choose {n/2}}2^{-n}
= \Theta(\sqrt{n})$ chips at position $\bv=0$.  The deterministic $P$-machine will give that
number with only constant error.  

As a further interpretation, consider a ``linear machine" in which a pile of $c$ chips are
split evenly, with $\frac{c}{2d}$ going in each direction.  The final number (which, in general,
will not be an integer) of chips at a position will simply be the expected number of chips in
the random walk machine.  We may consider the linear machine a ``relaxation'' of the $P$-machine,
 and our result states that this relaxation produces only constant ``error."
\section{The Proof}

We introduce some notation.  Chip distributions will be written as functions $\chi$ from $\Z^d$ to the 
nonnegative integers.  We let $n$ represent the total time of the process and let $t$ represent
the {\em time remaining} in the process so that $t$ counts backwards from $n$ to $0$.
Let $\chi_t$ denote the configuration of the $P$-machine with time remaining $t$.
Thus $\chi_n$ is the initial configuration and $\chi_0$ is the final
configuration.
We prescribe an ordering $\be_0,\ldots,\be_{2d-1}$ on 
those vectors (directions) with all coordinates zero except for one 
which is $+1$ or $-1$.
We fix
a map $\sigma : \Z^d \rightarrow \Z_{2d}$ denoting the initial settings of the rotors.  
Define $S_t$, a map from $\Z^d$ to $\R$, 
by setting $S_t(\bv)$ to be the probability that a random walk of time $t$
beginning at the origin ends at $\bv$.
For any distribution $\chi$ and any $t$ we define
\begin{equation}\label{Pnstar}
W_t(\bv,\chi) = \sum_{\bw\in Z^d}\chi(\bw)S_t(\bv-\bw),
\end{equation}
This $W_t$ acts as a weighted sum of all the chips when $t$ is the time remaining.  
Thus $W_0(\bv,\chi_0)$ is the actual number of chips at $\bv$ at the conclusion
of the $P$-machine process whereas $W_n(\bv,\chi_n)$ is the expected number
of chips at $\bv$ at the conclusion of the random walk process.

We (following \cite{L91}) write $\bv\lra n$ if the sum of the coordinates 
of $\bv$ has the same parity as $n$, and $\bw \lra \bw$ if $\bv$ and $\bw$ have the same parity in this sense.  Call a configuration $\chi$  ``even'' (``odd'') if, for all $\bv$ in the support of $\chi$, $\bv \lra 0$ (resp., $\bv \lra 1$).

Our main result is the following. 

\begin{theorem}\label{main} 
There is a constant $C_d$, depending on the dimension $d$ but not 
the initial even (or odd) configuration $\chi_n$, the total time $n$, nor the initial
rotor settings $\sigma$ so that
$$
 |W_0(\bv,\chi_0) - W_n(\bv,\chi_n)| < C_d.
$$
\end{theorem}

\begin{proof}
Clearly, we may assume without loss of generality that the configuration is even.  Then, for $n\geq t\geq 0$ we define
\begin{equation}\label{martingale}
 X_t = W_{t}(\bv,\chi_t)  
\end{equation}

We may interpret $X_t$ as the
expected number of chips at $\bv$ at the conclusion of a process in which
the first $n-t$ steps are
by the deterministic $P$-machine and the remaining $t$ steps are in the random
walk model.  The sequence $X_t$, $n\geq t\geq 0$, is then reminiscent of a martingale.
We first expand ``in time":
\begin{equation}\label{timeexpansion}
X_n - X_0 = \sum_{t=1}^n X_{t}-X_{t-1} = \sum_{t=1}^n W_{t}(\cdot,\chi_{t})-W_{t-1}
(\cdot,\chi_{t-1})
\end{equation}

The summand in (\ref{timeexpansion}) is the change in weight during 
one step of the $P$-machine.  
To analyze this change we set
\begin{equation} \label{changeweight} T^j(\bw,t) = S_t(\bw) - S_{t-1}(\bw - \be_j)  \end{equation}
This represents the change in weight of a single chip that was at position $\bv+\bw$ with
time remaining $t$ and with the rotor at that position set to state $j$.  Observe, critically,
$S_t(\bw)$ is the average of all possible $S_{t-1}(\bw-\be_j)$ so that 
\begin{equation} \label{harmonic}
\sum_{j=0}^{2d-1} T^j(\bw,u)  = 0  \end{equation}
Let $\Del_{\bw}$ denote the total contribution to $X_n-X_0$ from all chips 
moved from position $\bv+\bw$.  (Effectively we are here moving the origin
to $\bv$.)
Suppose that during the process a chip was moved from $\bv+\bw$ a 
total of $M$ times
and that the $i$-th time occured (starting the count at zero) when the remaining time was $t_i$.  Let
$c$ denote the initial state of position $\bv+\bw$.  Then
\begin{equation} \label{sumDelta} 
\Del_{\bw} = \sum_{i=0}^{M-1} T^{c+i}(\bw,t_i)  \end{equation}
(As before, we consider state $c+i$ as an element of $\Z_{2d}$.)
Reversing the sum of (\ref{timeexpansion}) we express
\begin{equation} \label{newsum}
X_n-X_0 = \sum_{\bw\in \Z^d}\Del_{\bw}.  \end{equation}

\subsection{Dimension One}

The case of $d=1$ is somewhat special, in that it requires separate arguments and is also susceptible to a more precise analysis.  Nonetheless, it closely resembles the higher dimensional cases, and so we produce the proof in the one-dimensional case first.

All vector quantities are simply integers when $d=1$, and we may write
$$
T^j(w,t) = S_t(w) - S_{t-1}(w \pm 1)
$$
where the choice of sign is $(-1)^j$.  Since, by (\ref{harmonic}), $T^0(w,t) = -T^1(w,t)$, $\Delta_w$ can be written as an alternating sum of $T^1(w,t_i)$ for some nonincreasing sequence $t_i$.

Now, we wish to show that the sequence $\{T^j(w,t)\}_{t \geq 1}$ is unimodal.  The unimodality of the sequence $T^1(w,t)$ is important because it implies that $T^1(w,t)$ is the concatenation of two monotone subsequences.  An alternating sum of any monotone sequence which has all the same sign is bounded by its largest term.  Therefore, we may estimate $\Delta_w$ by twice the maximum of $T^1(w,t)$.

It clearly suffices to consider the cases of $w \geq 1$ and $w=0$ -- otherwise write $\Delta_w$ in terms of $T^0$ instead of $T^1$.  Then, for $t \lra w$ with $w \geq 1$,
\begin{align}
T^1(w,t) &= 2^{-t} \binom{t}{(t+w)/2} - 2^{-t+1} \binom{t-1}{(t+w)/2-1}  \nonumber \\
& = 2^{-t+1} \binom{t-1}{(t+w)/2-1} \left (\frac{-w}{t+w} \right), \label{simpleformula}
\end{align}
which is always negative.  To determine when $T^j(w,t)$ is increasing or decreasing, we consider the ratio of consecutive terms in the range $t \geq w$:
\begin{align*}
\frac{T^1(w,t+2)}{T^1(w,t)} &= \frac{t+w}{4(t+w+2)} \cdot \frac{4t(t+1)}{(t+w)(t-w+2)} \\
& = \frac{t(t+1)}{(t+2)^2-w^2}.
\end{align*}
When $t(t+1) \leq (t+2)^2 - w^2$, i.e., when $t \geq (w^2-4)/3$, the quantity $T^1(w,t)$ is increasing; when $t \leq (w^2-4)/3$ it is decreasing.  (For $t < w$, the sequence is identically $0$.)  Hence, $T^1(w,t)$ is unimodal for those $t$ with $t \lra w$.

To compute the maximum of $T^1(w,t)$, we plug in $t = w^2/3 + O(1)$, getting $| T^1(w,t) | = O(w^{-2})$.  So, any alternating sum of $T^1(w,t_i)$ for $t_i$ nonincreasing satisfies this same bound.  (It is easy to see that, in the case of $w=0$, $T^1(w,t) = 2^{-t-1} \binom{t}{t/2}$ is also unimodal and bounded by a constant.)  Therefore, $\sum_{w \in \Z} O(w^{-2})$ converges, and Theorem \ref{main} follows in the case of $d=1$, by (\ref{newsum}).

\subsection{Higher Dimensions}

In dimensions $d \geq 2$, the situation is more complicated.  For $d\geq 3$
such elegant formulas as (\ref{simpleformula}) simply do not exist, and 
we must use estimates for the terms in question.  We suppress an argument for 
$d=2$ using exact formulae as our methods below apply to all $d\geq 2$.

\begin{claim} \label{claimone}
$\Del_{\bw} = O(|\bw|^{-(d+1)}\ln^{d-1}|\bw|)$.  That is, there exists a constant $K$ (depending
only on the dimension) such that $|\Del_{\bw}| \leq K|\bw|^{-(d+1)}\ln^{d-1}|\bw|$ for any choice of sequence 
$t_0,\ldots,t_{M-1}$ and initial state $c$ and any 
$\bw\in \Z^d$ with sufficiently high norm.  Further, for any fixed $\bw$ there is
a constant bound on $|\Del_{\bw}|$, independent of the choice of sequence $t_0,\ldots,t_{M-1}$.
\end{claim}

The case when $\bw$ is fixed is relatively simple and is given by Claim \ref{claimsix}.\\

Consider the expansion (\ref{sumDelta}) of $\Del_{\bw}$.  For each $i$ with $c+i=0$, apply
(\ref{harmonic}) and replace $T^0(\bw,t_i)$ with $\sum_{j\neq 0} -T^j(\bw,t_i)$.  Now for
each state $j\neq 0$ consider the sum of all terms with $T^j$.  When $c+i=j$ these will
be $T^j(\bw,t_i)$ and when $c+i=0$ these will be $-T^j(\bw,t_i)$.   These signs will
alternate, so that the sum can be written
\begin{equation} \label{alternating}
T^j(\bw,t_1)-T^j(\bw,t_2)+\ldots \pm T^j(\bw,t_s)
\end{equation}
where the $t_1,t_2,\ldots$ form a nonincreasing sequence.  (The initial sign might be negative, but that will not matter.)  Removing copies of the same
term we can assume the $t_1,t_2,\ldots$ form a decreasing sequence.  
(In particular, suppose the number of chips on a given position at a given time is
divisible by $2d$.  Then for each state $j$ the positive and negative contributions
will cancel out.  This is to be expected, as when the chips are distributed evenly
the weight function will not change.)
We have expressed
$\Del_{\bw}$ as the sum of a constant ($2d-1$) of these sequences so that to show
Claim \ref{claimone} it suffices to show

\begin{claim} \label{claimtwo}
$ T^j(\bw,t_1)-T^j(\bw,t_2)+\ldots \pm T^j(\bw,t_s) = O(|\bw|^{-(d+1)}\ln^{d-1}|\bw|) $  and
is bounded by a constant for any fixed $\bw$.
\end{claim}

As a further reduction we expand 
$$ T^j(\bw,t) = S_t(\bw) - S_{t-1}(\bw - \be_j)  = \frac{1}{2d}\sum_k S_{t-1}(\bw-\be_k)
- S_{t-1}(\bw-\be_j)  $$
where $k$ ranges over the $2d$ states.  The sum of Claim \ref{claimtwo} splits into
$2d-1$ sums. (When $k=j$ we get zero.) For fixed distinct states $j,k$ we set
$\bx = \bw-\be_j$ and $\bz = \be_j-\be_k$ and, again following  \cite{L91},
$$  \nabla_{\bz}(\bx,t) = S_t(\bx+\bz)-S_t(\bx) $$
As $|\bx|,|\bw|$ differ by at most one they are asymptotically the same and so
Claim \ref{claimtwo} reduces to showing

\begin{claim} \label{claimthree}
$ \nabla_{\bz}(\bx,t_1) - \nabla_{\bz}(\bx,t_2)+\ldots \pm \nabla_{\bz}(\bx,t_s) =
O(|\bx|^{-(d+1)}\ln^{d-1}|\bx|) $  and is bounded by a constant for any $\bx$.
\end{claim}
This we shall show for any fixed $\bz\lra 0$ and any decreasing sequence $t_1,\ldots,t_s$.

First, we show that the contribution from small $t$ (i.e., near the conclusion of
the process) is small.  With foresight
we set
\begin{equation}\label{defsmall} T = 
 K|\bx|^2/\ln^2(|\bx|) \end{equation}
where $K$ is a large constant to be specified.  $T$ will be our cutoff between
small and not small.
A result from \cite{L91} will be useful here: 
For some constant $c>0$, a random walk beginning at the origin ends at a point at least a 
distance $\alpha t^{1/2}$ away at time $t$ with probability at most $c e^{-\alpha}$.  
Now, $|\nabla_{\bz}(\bx,t)|$ is certainly bounded by the probability that the walk ends at least a 
distance $|\bx|$ from the origin plus the probability that it ends at least $|\bx+\bz|$ 
from the origin, so
\begin{align*}
|\nabla_{\bz}(\bx,t)| & \leq ce^{-|\bx|/t^{1/2}} + ce^{-|\bx+\bz|/t^{1/2}} \\
& \leq 2ce^{-(|\bx|-|\bz|)/t^{1/2}}.
\end{align*}
Now consider the the contribution from the small $t<T$. 
We select $K$ so that all $|\nabla_{\bz}(\bx,t)| = O(|\bx|^{-d-3})$.  The total
contribution to the sum of Claim \ref{claimthree} from all $t\leq T$, there being
fewer than $|\bx|^2$ terms, is $O(|\bx|^{-d-1})$.  

We can now assume that all $t_i \geq T$.
We approximate the $\nabla_{\bz}(\bx,t)$ by the following special case of the ``Local
Central Limit Theorem" for simple random walks (q.v. \cite{L91}):

\begin{equation} \label{LCLT}
| \nabla_{\bz}(\bx,t) - (p(t,\bx+\bz) - p(t,\bx)) | = |\bx|^{-2} O(t^{-(d+1)/2})
\end{equation}
where $|\cdot|$ denotes the $L^2$-norm, $t\lra \bx$,
and 
\begin{equation}\label{defp}
p(t,\bx):= 2 (\frac{d}{2 \pi t})^{d/2} e^{-d|\bx|^2/2t} \end{equation}
(Here, and throughout the rest of the proof, all hidden constants may depend on $\bz$.  This
does not affect our original problem since $\bz$ takes on only a finite set of values.)

An alternating sum of $\nabla_{\bz}(\bx,t)$ is then, in absolute value, at most the
absolute value of the alternating sum of the approximation $p(t,\bx+\bz)-p(t,\bx)$ and
the sum of the absolute values of the errors.   The sum
of the absolute values of the errors for $t\geq T$ is 
then $O(|\bx|^{-2}T^{(-d+1)/2})$.  (Note that this argument would {\em not}
work for $d=1$ as the bounds on the errors would form a divergent harmonic
series.) 
Our choice of cutoff $T$ ensures that this is $O(|\bx|^{-(d+1)}\ln^{d-1}|\bx|)$.

Now, we wish to show that an alternating sum of the 
quantity $f(t) = p(t,\bx) - p(t,\bx+\bz)$ is $O(|\bx|^{-(d+1)})$.
We regard this as the ``main term" of our calculations.  This is
the critical place where we need the alternation of the sum.
That the sum alternates is a direct consequence of the 
construction of the $P$-machine, and represents in our minds
the ``cancelling out process" when a single rotor distributes
chips in a (relatively) even manner.

\begin{claim}\label{claimfour} $f(t)= O(|\bx|^{-(d + 1)})$ \end{claim}

Let $\rho = |\bx+\bz| - |\bw|$, so $|\rho| \leq |\bz|$.  Since $f(t)$ is zero when $\min\{|\bx|,|\bx+\bz|\} > t$, we may write
\begin{align*}
f(t) & \ll t^{-d/2} | e^{-d|\bx|^2/2t} - e^{-d|\bx+\bz|^2/2t} | \\
&= t^{-d/2} e^{-d|\bx|^2/2t} |1 - e^{-d \rho (|\bx+\bz| + |\bx|)/2t}| \\ 
& \ll t^{-d/2} e^{-d|\bx|^2/2t} d|\rho| (|\bx+\bz|+|\bx|)/{2t} \\
& \ll t^{-d/2 - 1} |\bx| e^{-d|\bx|^2/2t}.
\end{align*}
It is a matter of elementary calculus to show that this function is 
maximized when $t = \frac{d}{d+2} |\bx|^2$, in which case $f(t)$ satisfies Claim 
\ref{claimfour}.

\begin{claim}\label{claimfive}  The sequence $f(t)$ has at most $6$ local extrema.
(We say that $t$ is a local extremum if $f(t)>\max[f(t+2),f(t-2)]$ or $f(t)<\min[f(t+2),f(t-2)]$).
\end{claim}

Letting $B = 2 (d/2 \pi)^{d/2}$,, we may write
$$
p(t,\bx)-p(t,\bx+\bz) = B t^{-d/2} (e^{-d|\bx|^2/2t} - e^{-d(|\bx+\bz|^2/2t}) \\
$$
Setting the derivative equal to zero gives
\begin{equation} \label{firstderivative}
e^{-\gamma y} = K_1 + \frac{K_2}{y-1}
\end{equation}
after making the substitutions $R = |\bx|^2$, $S = |\bx+\bz|^2$, $y = R/t$, $K_1 = 
S/R$, $K_2 = (S-R)/R$, and $\gamma = (R-S)/2dR$.  This splits into the 
two regions $y < 1$ and $y > 1$. In each region the number of zeroes is at most 
one more than the number of zeroes of the derivative.  Therefore, setting 
the derivative of (\ref{firstderivative}) equal to zero gives the equation
$$
e^{-\gamma y} (y-1)^2 - K_3 = 0
$$
The number of zeroes of this is at most (in each region) the number of zeroes of 
its derivative, which is
$$
e^{-\gamma y} [2(y - 1) - \gamma (y - 1)^2] = 0.
$$
This equation clearly has at most two solutions, and the proof of Claim 
\ref{claimfive} is complete.

From Claim \ref{claimfive} the sequence $f(t)$ can be split into a bounded
$(\leq 7)$ number of monotone sequences.  {\em On each such sequence an alternating
sum is at most twice the maximal absolute value of the terms.}  Applying
Claim \ref{claimfour} the entire sum is $O(|\bx|^{-(d+1)})$ which concludes
Claim \ref{claimthree} and thus Claim \ref{claimone} in the asymptotic case
when $\bx\ra\infty$.

\begin{claim}\label{claimsix} For any fixed $\bw$ there is
a constant bound on $|\Del_{\bw}|$.
\end{claim}
The reductions to Claim \ref{claimthree} are as given.   The sequence
$p(t,\bx+\bz)-p(t,\bx)$ clearly has absolute value at most two and the
argument of Claim \ref{claimfive} holds for any $\bx$ so an alternating
sum of the $p(t,\bx+\bz)-p(t,\bx)$ is bounded.  The error terms, from
(\ref{LCLT}), are $O(t^{-(d+1)/2})$ and so their sum is bounded for all $d\geq 2$.

Having shown Claim \ref{claimone}, we may now show Theorem \ref{main}.
Equation (\ref{newsum}) and this Claim gives 
\begin{equation}\label{laststep}
X_n-X_0=O(\sum_{\bw\in Z^d}|\bw|^{-(d+1)}\ln^{d-1}|\bw|) = O(1).
\end{equation}

\end{proof}

\section{Conjectures}

In the course of proving the main theorem, we encountered the following 
intriguing conjecture.

\begin{conjecture} \label{conjectureone} For each $\xi : \Z^d \rightarrow \R$ with finite support, and for all $\bv \in \Z^d$, the 
function $W_{2n}(\bv,\xi)$ has at most $K$ sign changes in the variable $n$, 
where $K$ depends on $\xi$ and $d$ but not on $\bv$.
\end{conjecture}

To see where this statement comes from, note that, for any $\xi$, we may always use the harmonic property to write $W_{2n+2}(\bv,\xi) - W_{2n}(\bv,\xi) = W_{2n}(\bv,\xi^\prime)$ for some $\xi^\prime : \Z^d \rightarrow \R$.  Then, $W_{2n}(\bv,\xi)$ has a bounded number of local extrema if and only if $W_{2n}(\bv,\xi^\prime)$ has a bounded number of sign changes.  Therefore, if this conjecture were true, we could simply estimate the alternating sums in the proof by a constant times their largest term, and the burden of bounding the error terms in (\ref{LCLT}) would be reduced somewhat.  Our intuition is that this statement is true because ``oscillations'' in $p(n,\bv)$ occur only to the modulus $2$.

Furthermore, we believe that the polylogarithmic factor in the bound on $|\Delta_{\bw}|$ is unnecessary.  Indeed, this would follow from Conjecture \ref{conjectureone}. 

\begin{conjecture} For any $d$, $|\Delta_\bw | = O(|\bw|^{-(d+1)})$.
\end{conjecture}

For very simple $\xi$, such as one whose support consists of a single vertex, we actually expect unimodality.  That is,

\begin{conjecture} \label{conjecturetwo} The probability that a random walk beginning at the origin in $Z^d$ 
arrives at $\bv$ at time $n \leftrightarrow \bv$ is a unimodal function of $n$.
\end{conjecture}

The estimates given above are inadequate to prove this simple statement.  Furthermore, attempting to 
show unimodality via log-concavity is futile, since these functions are often 
(perhaps always?) log-convex in the tail.  In the $d=1$ case, we can show unimodality by brute force.  Indeed, $p(n,x) = 2^{-n} \binom{n}{(n+x)/2}$ for $x \geq 0$, and the ratio of consecutive terms is given by
$$
\frac{p(n,x)}{p(n-2,x)} = \frac{n(n-1)}{n^2 - x^2}
$$
so the sequence is $0$ for $n < x$, increasing up to $n = x^2$, and decreasing thereafter.  A similar argument also works for $d = 2$.  It is clear that it suffices to prove unimodality for the half-quadrant where $x \geq 0$ and $y \leq x$.  Note that, in a random walk on $\Z^2$, the quantities $x+y$ and $x-y$ vary {\it independently} at each step.  Furthermore, $x+y$ and $x-y$ follow one-dimensional simple random walks.  Therefore, we may decompose the two-dimensional walk as follows:
\begin{align*}
p(n,(x,y)) &= p(n,x+y) \cdot p(n,x-y) \\
& = 4^{-n} \binom{n}{(n+x+y)/2} \binom{n}{(n+x-y)/2}.
\end{align*}
If we consider the ratio of successive terms,
$$
\frac{p(n,(x,y))}{p(n-2,(x,y))} = \frac{n^2(n-1)^2}{(n+x+y)(n+x-y)(n-x+y)(n-x-y)}.
$$
In the range $n > x+y$, this quantity is $\leq 1$ precisely when $g(n) = 2n^3 - rn^2 + s$ is positive, where $r = 2x^2 + 2y^2 + 1$ and $s=(x^2-y^2)^2$.  When $n = x+y$,
$$
g(x+y) = - (x+y)^2 (x+y-1)^2 \leq 0.
$$
On the other hand, $g^\prime$ has two roots: one at $n=0$ and one at $n = r/3$.  Since $0 \leq x+y$, at most one root of $g$ occurs at in $[x+y,+\infty)$, so we may conclude that $p(n,(x,y))$ is unimodal in $n$.

Unfortunately, this kind of analysis does not seem to be possible in higher dimensions, so further insights are necessary to resolve the conjecture when $d \geq 3$.

\section{Acknowledgements}

Thank you to Noam Berger, Jim Propp, and David Wilson for helpful discussions and 
insights.

\end{document}